\documentclass[11pt,a4paper]{article}
\usepackage{amsmath,amsfonts,amssymb,mathrsfs}
\usepackage{graphicx}
\usepackage{subfigure}
\usepackage{pstricks}
\usepackage{pst-node}
\usepackage{xcolor}
\usepackage{url}
\newrgbcolor{blu}{.1 0 .7}
\newrgbcolor{lgrey}{.9 .9 .9}
\newrgbcolor{grey}{.7 .7 .7}
\newrgbcolor{dgrey}{.3 .3 .3}
\newrgbcolor{lred}{1 .8 .8}
\newrgbcolor{mlred}{1 .7 .7}
\newrgbcolor{bblue}{.6 .6 1}
\newrgbcolor{lblue}{.8 .8 1}
\newrgbcolor{mblue}{.2 .25 1}
\newrgbcolor{dblue}{.3 .32 1}
\newrgbcolor{ddblue}{.4 .4 1}
\newrgbcolor{lgreen}{.8 1 .8}
\newrgbcolor{mgreen}{ .7 1 .7}
\newrgbcolor{ggreen}{ .1 .75 .35}
\newtheorem{theorem}{\bf Theorem}[section]

\newtheorem{corollary}[theorem]{\bf Corollary}
\newtheorem{definition}[theorem]{\bf Definition}
\newtheorem{remark}[theorem]{\bf Remark}
\newtheorem{proposition}[theorem]{\bf Proposition}

\newcommand{\R}{\mathbb{R}}

\def \rn {{\mathbb {R}}^{N}}
\def \rnn {{\mathbb {R}}^{N+1}}

\def \B {{\cal B}}
\def \K  {K}  
\def \L {\mathscr{L}}
\def \F {\mathscr{F}}
\def \A {\mathscr{A}}
\def \P {{\mathscr{P}}}

\def\p{\partial}

\def \o {{\omega}}

\def \g {{\gamma}}
\def \d {{\delta}}

\def \eps {{\varepsilon}}
\def \epsilon {{\varepsilon}}

\def \l {{\lambda}}

\def \t {{\tau}}

\def \phi {{\varphi}}
\def \G {{\Gamma}}
\def \O {{\Omega}}

\usepackage[english]{babel}
\begin{document}
\title{Harnack inequalities and Bounds for Densities of Stochastic Processes}
\author{{\sc{Gennaro Cibelli}
\thanks{Dipartimento di Scienze Fisiche, Informatiche e Matematiche, Universit\`{a} di Modena e Reggio Emilia, Via
Campi 213/b, 41125 Modena (Italy). E-mail: gennaro.cibelli@unimore.it}
\sc{Sergio Polidoro}
\thanks{Dipartimento di Scienze Fisiche, Informatiche e Matematiche, Universit\`{a} di Modena e Reggio Emilia, Via
Campi 213/b, 41125 Modena (Italy). E-mail: sergio.polidoro@unimore.it}
}}
\maketitle

\footnotesize

{\it \hfill Dedicated to Valentin Konakov in occasion of his 70th birthday.}

\bigskip

\begin{abstract}
We consider possibly degenerate parabolic operators in the form
$$
  \L  = \sum_{k=1}^{m}X_{k}^{2}+X_{0}-\p_{t},
$$
that are naturally associated to a suitable family of stochastic differential equations, and satisfying the
H\"ormander condition.
Note that, under this assumption, the operators in the form $\L$ have a smooth fundamental solution that agrees with the
density of the corresponding stochastic process. We describe a method based on Harnack inequalities and on the
construction of Harnack chains to prove lower bounds for the fundamental solution. We also briefly discuss PDE and SDE
methods to prove analogous upper bounds. We eventually give a list of meaningful examples of operators to which the
method applies
\end{abstract}

\normalsize\

\section{Introduction}\label{intro}

Let  $(W_t)_{t\geq 0}$ denote an $m$-dimensional Brownian motion, $W_t=(W_t^1,\ldots, W_t^m)$ on some filtered probability space $(\O, \F, (\F)_{t \ge 0}, \mathbb{P})$. We consider a collection of space-time functions $(\sigma_{ij})_{(i,j) \in \{1,...,N\} \times \{1,...,m\}}, \ (b_i)_{i \{1,...,N\}}$ such that the following SDE
\begin{equation}\label{diffusion}
 dZ^i_t=\sum_{j=1}^m\sigma_{ij}(Z_t,t) \circ dW^j_t+ b_i(Z_t,t)dt, \qquad i= 1, \dots, N, \quad t \geq 0
\end{equation}
is well posed at least in the weak sense. Here ``$\circ \, dW_t$'' stands for the Stratonovich integral. We denote by $Z^{x_0}_t$
the solution of the SDE (\ref{diffusion}) with initial condition $Z_0^{x_0}={x_0}$. The equation \eqref{diffusion} is
associated to the Kolmogorov operator
\begin{equation*} 
  \L= \sum_{i=1}^mX_i^2 + X_0-\partial_t,
\end{equation*}
where
\begin{equation}\label{campi}
  X_i(x,t)=\tfrac{1}{\sqrt 2}\sum_{j=1}^m \sigma_{ij}(x,t)\p_{x_j}, \ i=1, \dots,m, \qquad
  X_0(x,t)= \sum_{i=1}^N b_i(x,t)\partial_{x_i}.
\end{equation}

In this note we describe a general method to prove upper and lower bounds for the fundamental solution of $\L$. Specifically, we say that a non-negative function $\Gamma(x,t; y,s)$ defined for
$x,y\in\R^{N}$ and $t>s$, is a fundamental solution for $\L$ if:
\begin{itemize}
  \item[i)] in the weak sense, $\L\G(\cdot,\cdot;y,s)=0$ in $]s,+\infty[\times \R^{N}$ and $\L^{*}\G(t,x;\cdot,\cdot)=0$  in $]-\infty,t[\times
  \R^{N}$ where $\L^*$ denotes the formal adjoint operator of $\L$;
  \item[ii)] for any bounded function $\phi\in C(\R^{N})$ and $x,y\in\R^{N}$, we have
\begin{align}
  \lim_{(x,t)\to(y,s)}u(x,t)=\phi(y), 
  \qquad \lim_{(y,s)\to(x,t)}v(y,s)=\phi(x), 
\end{align}
where
  \begin{align}\label{ae11}
 u(x,t):=\int_{\R^{N}}\G(x,t;y,s)\phi(y)dy,\qquad 
 v(y,s):=\int_{\R^{N}}\G(x,t;y,s)\phi(x)dx.
\end{align}
\end{itemize}

Note that the functions in \eqref{ae11} are weak solutions of the following backward and forward Cauchy
problems:
  $$
  \begin{cases}
    \L u(t,x)= 0,\ &(x,t)\in\, ]s, +\infty[\times \R^{N}, \\
    u(x,s)=\phi(x),&x\in\R^{N},
  \end{cases}\quad
  \begin{cases}
    \L^{*}v(y,s)= 0,\ & (y,s)\in\,]-\infty,t[\times \R^{N}, \\
    v(y,t)=\phi(y), & y\in\R^{N}.
  \end{cases}
  $$

We introduce the $N \times N$ matrix $A(x,t)=\left(a_{ij}(x,t)\right)_{i,j = 1 ,\dots, N}$ whose elements are
$$
  a_{ij}(x,t)=\tfrac{1}{2}\sum_{k=1}^m \sigma_{ik}(x,t)\sigma_{jk}(x,t), \qquad i,j= 1, \dots, N,
$$
and we note that
\begin{equation*}
  \left\langle A(x,t) \xi, \xi \right \rangle = \tfrac{1}{2} \| \sigma(t,x) \xi  \|^2 \ge 0.
  \quad \text{for every} \quad \xi \in \R^N.
\end{equation*}
If the smallest eigenvalue of $A(t,x)$ is uniformly positive we say that the operator $\L$ is uniformly
parabolic.

 A keystone result in the theory of parabolic partial differential equations reads as follows. Assume that
there exist two positive constants $\lambda, \Lambda$ such that
\begin{equation} \label{e-unif-pos}
  \lambda |\xi|^2 \le
  \left\langle A(x,t) \xi, \xi \right \rangle
  \le \Lambda |\xi|^2, \quad
  \text{for every} \ (x,t)  \in \R^N\times ]0,T[ , \ \text{and} \ \xi \in \R^N.
\end{equation}
If $\Gamma = \Gamma(x,t,\xi, \tau)$ denotes the fundamental solution of the PDE
\begin{equation} \label{e-unif-par}
  \partial_{t} u (x,t) = \sum_{i,j=1}^N \partial_{x_i} \left( a_{ij}(x,t) \partial_{x_j} u(x,t) \right), \qquad (x,t)
\in \R^N \times ]0,T[,
\end{equation}
then there exist positive constants $c^-, C^-, c^+,C^+$ only depending on $N, \Lambda, \lambda$ such that
\begin{equation} \label{e-keystone}
  \frac{c^-}{(t-\tau)^{N/2}} \exp \left(- C^- \frac{|x-\xi|^2}{t-\tau} \right) \le \Gamma(x,t;\xi,\tau) \le
  \frac{C^+}{(t-\tau)^{N/2}} \exp \left(- c^+ \frac{|x-\xi|^2}{t-\tau} \right),
\end{equation}
for every $(x,t), (\xi,\tau) \in \R^N \times ]0,T[$ with $\tau<t$. We  emphasize that the constants in \eqref{e-keystone} do not depend on $T$. This upper bound has been proved by Aronson
\cite{Aronson} for operators with bounded measurable coefficients $a_{ij}$'s, while the lower bound has been
proved by Moser \cite{Moser3, Moser3bis}. The results by Aronson and by Moser improve the earliest estimates given by Nash in his seminal work \cite{Nash}. We also refer to the article of Krylov and Safonov \cite{KrylovSafonov2} for non-divergence form operators.

\smallskip

The results described above have been extended by several authors to possibly degenerate operators in the form
\begin{equation}\label{e1}
\L := \sum_{k=1}^{m}X_{k}^{2}+Y,  \qquad Y :=X_0-\p_t,
\end{equation}
where $X_{0}, X_{1}, \dots, X_{m}$ are smooth vector fields on $\R^{N+1}$, that is
\begin{equation}\label{smooth_vector}
  X_{i}(x,t)=\sum_{j=1}^{N}c_{i,j}(x,t)\p_{x_{j}}, \qquad i=0,\dots,m.
\end{equation}
for some smooth functions $c_{i,j}$'s. In particular, upper bounds have been proved by a PDE approach that goes back
to Aronson's work \cite{Aronson}, or by an approach based on Lyapunov functions (see \cite{MetafunePallaraRhandi}
and the references therein). Several authors prove bounds analogous to \eqref{e-keystone} in the framework of
stochastic processes. We refer to the works of Malliavin \cite{Malliavin}, Kusuoka and Stroock \cite{KusuokaStroock},
where a general method to prove upper bounds for density is introduced and to the work of Ben Arous and L\'eandre
\cite{leandre}, where the Malliavin Calculus is further developed. We also refer to the monograph of Nualart
\cite{Nualart} for a comprehensive presentation of this subject.

\smallskip

In general, lower bounds have been proved by following the idea introduced by Moser in \cite{Moser3}. In this note we
briefly describe this method for uniformly parabolic partial differential equations, then we give an overview of more
recent articles where it has been adapted to the study of degenerate parabolic equations in the form \eqref{e1}.
This idea is also used in the works where lower bounds are proved by probabilistic methods. We refer to Kohatsu-Higa
\cite{Kohatsu1}, Bally \cite{Bally2}, Bally and Kohatsu-Higa \cite{BallyHiga}.

\medskip

We now give a list of examples of operators considered in this note. Each one of them is the prototype of a wide family
of differential operators.
\begin{itemize}
\item \emph{Heat operator on the Heisenberg group} $\L =  X_1^2 +X_2^2  - \p_t,$ where
$$
  X_1=\p_x - \tfrac{1}{2}y \p_w, \qquad X_2=\p_y + \tfrac{1}{2}x \p_w.
$$
Note that $\L$  acts on the variable $(x,y,w,t) \in \R^{4}$, and writes in the form \eqref{e1} with $X_0=0$. The
degenerate elliptic operator $\Delta_{\mathbb{H}}=X_1^2+X_2^2$ is said \emph{sub-Laplacian} on the Heisenberg group.
\end{itemize}
\begin{itemize}
\item \emph{Kolmogorov Operator} $\L=\p_{xx} +x\p_y -\p_t, (x,y,t) \in \R^3$. In this case $\L = X^2 + Y$ with
$X=\p_x, \ Y=x\p_y-\p_t$.
\end{itemize}
\begin{itemize}
\item \emph{More Degenerate Kolmogorov Operators} $\L=\p_{xx} +x^2\p_y -\p_t, (x,y,t) \in \R^3$. In this
case $\L = X^2 + Y$ with $X=\p_x, \ Y=x^2\p_y-\p_t$.
\end{itemize}
\begin{itemize}
\item \emph{Asian Option Operator} $\L=x^2 \p_{xx} +x\p_x +x\p_y -\p_t, (x,y,t) \in \R^+ \times \R^2$. In
this case $\L = X^2 + Y$ with $X=x\p_x, \ Y=x\p_y-\p_t$.
\end{itemize}
All the operators in the above list are strongly degenerate, since the smallest eigenvalue of the characteristic form is
zero for all the above examples. In general, operators in the form \eqref{e1} cannot be uniformly parabolic if $m < N$.
On the other hand, all the examples do satisfy the following condition:

\smallskip

\noindent {\sc Hypothesis [H]}\ {\it $\L = \sum_{k=1}^{m}X_{k}^{2}+Y$ satisfies the H\"{o}rmander
condition if
\begin{equation*} 
  \emph{rank} \ (\text{\rm Lie}\{X_{1},\dots,X_{m},Y\}(x,t))=N+1, \qquad \text{for every} \quad (x,t) \in \rnn.
\end{equation*}}
In the sequel we only consider operators $\L$ satisfying the H\"ormander condition. It is know that, for this family of
operators, the law of the stochastic process \eqref{diffusion} is absolutely continuous with respect to the Lebesgue
measure in $\rn$, and that its density is smooth. Moreover, for every pairs $(\xi,\tau),\ (x,t) \in \R^N \times [0,T[$ with $\tau>t$, the density $p(\xi, \tau; x,t)$ is linked with the
fundamental
solution $\G$ of $\L$. Precisely, if $p$ denotes the density of the process
$$
\left\{
  \begin{array}{ll}
    dZ_s^i=\sum_{i,j=1}^m \sigma_{ij}(Z_s,T-s) \circ dW_s^i+b_i(Z_s,T-s)ds, \ \  &\hbox{$i=1,...,N, \quad t< s \le T$;} \\
    Z_t^i=x_i, & \hbox{$i=1,...,N$,}
  \end{array}
\right.
$$
then $\G$ is defined by the relation
$$
  \G(x,t;\xi,\tau)=p(\xi, T-\tau; x,T-t).
$$
It is known that the regularity properties of the operators satisfying the H\"ormander condition are related to a
Lie group structure that replaces the usual Euclidean one. In the proof of the lower bounds for positive solutions
the geometric aspects of this \emph{non Euclidean} structure will be explicitly used. To make  the exposition clear, in
Section \ref{sec2} we recall the method used by Moser in \cite{Moser3} to prove the lower bound in \eqref{e-keystone}
for uniformly parabolic operators. In Section \ref{sec3} we describe how the method outlined in Section \ref{sec2} is
adapted to the degenerate ones, satisfying the H\"ormander condition [H]. The remaining Sections \ref{sec4},
\ref{sec5}, \ref{sec6}  and \ref{sec7} are devoted to the examples listed above.

\section{Uniformly parabolic equations}\label{sec2}
In this section we describe the method introduced by Moser \cite{Moser3} to prove the lower bound \eqref{e-keystone} of
the fundamental solution for uniformly parabolic equations. The main ingredient of the method is the \emph{parabolic
Harnack inequality}, first proved by Hadamard \cite{Hadamard} and, independently, by Pini \cite{Pini} in 1954 for the
heat equation, then by Moser \cite{Moser3, Moser3bis} for uniformly parabolic equations in divergence form
\eqref{e-unif-par}. Its statement requires some notation (see Fig. 1). Let
\begin{equation*}
  Q_r(x,t) = B(x,r) \times ]t- r^2, t[ ,
\end{equation*}
denote the  parabolic cylinder whose upper basis is centered at $(x,t)$. Let $\alpha, \beta, \gamma, \delta \in ]0,1[$
be given constants, with $\alpha < \beta < \gamma <1$,
\begin{equation*}
  Q^-_r(x,t) = B(x,\delta r) \times ]t- \gamma r^2, t - \beta r^2[ \qquad
  Q^+_r(x,t) = B(x,\delta r) \times ]t- \alpha r^2, t[.
\end{equation*}
\begin{center}
\vskip -4mm
\begin{pspicture*}(-.9,.8)(6.1,4.1)
\scalebox{.80}{%
    \pcline[linecolor=dgrey, linestyle=dashed](3,3.2)(3,2.4)%
    \psframe[fillstyle=solid, fillcolor=lgrey, linecolor=dgrey](1.1,1.6)(4.9,2.4)%
    \pcline[linecolor=dgrey, linestyle=dashed](3,1.6)(3,2.39)%
    \psframe[fillstyle=solid, fillcolor=lgrey, linecolor=dgrey](1.1,3.2)(4.9,4.01)%
    \pcline[linecolor=dgrey, linestyle=dashed](3,3.2)(3,4)%
    \pcline[linecolor=black](0,4)(6,4)%
    \pcline[linecolor=black](0,1)(6,1)%
    \pcline[linecolor=black](6,4)(6,1)%
    \pcline[linecolor=black](0,4)(0,1)%
    \dotnode(3,4){O}\nput{90}{O}{$(x,t)$}
    \uput[0](6,2.6){$r^2$}
    \uput[0](3,2.8){{$(\beta - \alpha) r^2$}}
    \uput[-90](3.8,1.6){{ $\delta r$}}
    \uput[0](3,2){{$(\gamma- \beta) r^2$}}
    \uput[90](3.8,4){{ $\delta r$}}
    \uput[0](3,3.6){{ $\alpha r^2$}}
  }
\end{pspicture*}
\end{center}
\begin{center}
 {\scriptsize \sc Fig. 1 - Parabolic Harnack inequality.}
\end{center}

\begin{theorem} [{\sc Parabolic Harnack inequality}] \label{th-Harnack-P}
Let $Q_r(x,t) \subset \R^{N+1}$, and let $\alpha, \beta, \gamma, \delta \in ]0,1[$ be  given constants, with $\alpha <
\beta < \gamma <1$. Then there exists $C = C(\alpha, \beta, \gamma, \delta, \lambda, \Lambda, N)$ such that
\begin{equation*}
  \sup_{{Q^-_r(x,t)}} u \le C \inf_{{Q^+_r(x,t)}} u
\end{equation*}
for every $u: Q_r(x,t) \to \R, u \ge 0,$ satisfying \eqref{e-unif-par}. Here $\lambda, \Lambda$ are the constants in
\eqref{e-unif-pos}.
\end{theorem}

\begin{remark} \label{rem-Harnack-inv}
Note that  $C$ does not depend on the point $(x,t)$ and on $r$, then the Harnack inequality is invariant with respect
to the \emph{Euclidean translation} $(x,t) \mapsto (x + x_0, t+ t_0)$, and to the \emph{parabolic dilation} $(x,t)
\mapsto (r x,r^2 t)$. For this reason, the above statement is often referred to as \emph{invariant Harnack inequality}.
\end{remark}

In the sequel we will use the following version of the parabolic Harnack inequality  (see Fig. 2). For any given $c \in
]0,1[$ we
denote by
\begin{equation*}
  P_r(x,t) =  \big\{(y,s) \in Q_r(x,t) \mid 0 < t-s \le  c r^2 < t, |y-x|^2 \le t-s \big\}.
\end{equation*}
\begin{center}
\begin{pspicture*}(-.9,.8)(6.1,4.1)
  \scalebox{.8}{%
    \parabola[linecolor=dgrey,fillstyle=solid,fillcolor=lgrey](1,2)(3,4)
    \pcline[linecolor=dgrey](1,2)(5,2)%
    \dotnode(3,4){O}\nput{90}{O}{$(x,t)$}
    \pcline[linecolor=black](0,4)(6,4)%
    \pcline[linecolor=black](0,1)(6,1)%
    \pcline[linecolor=black](6,4)(6,1)%
    \pcline[linecolor=black](0,4)(0,1)%
}
\end{pspicture*}
\end{center}
\begin{center}
  {\scriptsize \sc Fig. 2 - Parabolic Harnack inequality.}
\end{center}

\begin{corollary}  \label{cor-Harnack-P}
Let $Q_r(x,t) \subset \R^{N+1}$, and let $c \in ]0,1[$ be  a given constant. Then there exists $C = C(c, \lambda,
\Lambda, N)$ such that
\begin{equation*}
  \sup_{{P_r(x,t)}} u \le C  u (x,t)
\end{equation*}
for every $u: Q_r(x,t) \to \R, u \ge 0,$ satisfying \eqref{e-unif-par}. Here $\lambda, \Lambda$ are the constants in
\eqref{e-unif-pos}.
\end{corollary}

\medskip

\noindent {\sc Proof.}  For every positive $\rho$ we denote
\begin{equation*}
  S_\rho(x,t) = B(x, \rho) \times \{t - \rho^2\}.
\end{equation*}
Let $\alpha, \beta, \gamma \in ]0,1[$ be such that $\alpha < \beta \le c \le \gamma <1$, and let $\delta = \sqrt{c}$.
Then, for every $\rho \in [0,r]$ we have that $u$ is a non-negative solution of \eqref{e-unif-par} in the
domain $Q_\rho(x,t)$. Since $S_\rho(x,t) \subset Q^-_\rho(x,t)$,  from Theorem \ref{th-Harnack-P} we obtain
\begin{equation*}
  \sup_{{S_\rho(x,t)}} u \le \sup_{{Q^-_\rho(x,t)}} u \le C \inf_{{Q^+_\rho(x,t)}} u \le C u(x,t),
\end{equation*}
and the conclusion follows from the fact that $P_r(x,t) = \cup_{0 < \rho \le r} S_\rho(x,t)$.
\hfill $\square$

\medskip

With Corollary \ref{cor-Harnack-P} in hand, we can easily obtain the following \emph{non local} Harnack inequality,
first proved by Moser  (Theorem 2 in \cite{Moser3}). We also refer to  Aronson \& Serrin \cite{AronsonSerrin} for more
general uniformly parabolic differential operators.

\begin{theorem} \label{th-Harnack-NL}
Let $u: \R^N \times ]0,T[ \to \R$ be a non-negative solution of the parabolic equation \eqref{e-unif-par}. Then there
exists a positive constant $C = C(c, \lambda, \Lambda, N)$ such that
\begin{equation*}
  u(x,t) \le C^{1 + \frac{|x_0-x|^2}{t_0-t}} u(x_0,t_0),
\end{equation*}
for every  $(x_0,t_0), (x,t) \in \R^N \times ]0,T[$ with $t_0 - t < c \, t_0$.
\end{theorem}

\medskip

\noindent {\sc Proof.} Let $(x_0,t_0), (x,t) \in \R^N \times ]0,T[$, with $t_0 - t < c t_0$, choose $r =
\sqrt{t}$ and note that the cylinder $Q_r(x_0,t_0)$ is contained in $\R^N \times ]0,T[$. If $(x,t) \in P_r(x_0,t_0)$ we
simply apply Corollary \ref{cor-Harnack-P} and the proof is complete. If otherwise $(x,t) \not \in P_r(x_0,t_0)$, we
consider the segment whose end points are $(x_0,t_0)$ and $(x,t)$, and denote by $(x_1, t_1)$ the point where it
intersects the boundary of $P_r(x_0,t_0)$. Note that $t_1 \ge  t > (1-c) t_0$, then  $(x_1, t_1)$ belongs to the
lateral part of the boundary of $P_r(x_0,t_0)$. By Corollary \ref{cor-Harnack-P} we have
\begin{equation*}
  u (x_1, t_1) \le C  u (x_0,t_0).
\end{equation*}
We then iterate the argument. We define a finite sequence $(x_j, t_j)$, with $j=2, \dots, k$ such that $(x_j, t_j)$
belonging to the boundary of $P_r(x_{j-1},t_{j-1})$ for  $j=2, \dots, k$, and $(x,t) \in P_r(x_k,t_k)$  (see Fig. 3).
By applying $k$ times Corollary \ref{cor-Harnack-P} we then find
\begin{equation*}
  u(x,t) \le C u (x_k, t_k) \le C^2 u (x_{k-1}, t_{k-1}) \le \dots \le C^{k + 1}  u (x_0,t_0).
\end{equation*}
\begin{center}
\begin{pspicture*}(-.9,0)(8.1,3.3)
  \scalebox{.7}{%
    \parabola[linecolor=dgrey,fillstyle=solid,fillcolor=lgrey](-1,0)(2,4)
    \dotnode(2,4){O}\nput{90}{O}{$(x_0,t_0)$}
    \dotnode(7,1){P}\nput{60}{P}{$(x,t)$}
    \dotnode(2,4){O}
    \dotnode(7,1){P}
    \dotnode(3.27,3.25){Z1}\nput{60}{Z1}{$(x_1, t_1)$}
    \multirput(0,0)(1.25,-.75){4}{\parabola[linecolor=dgrey,fillstyle=solid,fillcolor=lgrey](-1,0)(2,4)}
    \pcline[linecolor=dgrey](2,4)(7,1)%
    \dotnode(4.5,2.5){Z2}\nput{60}{Z2}{$(x_2, t_2)$}
    \dotnode(5.75,1.75){Z3}\nput{60}{Z3}{$(x_3, t_3)$}
    \dotnode(3.27,3.25){Z1}
    \dotnode(2,4){O}
    \dotnode(7,1){P}
}
\end{pspicture*}
\end{center}
\begin{center}
  {\scriptsize \sc Fig. 3 - Harncak chain.}
\end{center}
To conclude the proof it is sufficient to note that the integer $k$ only depends on the slope of the line connecting
$(x_0,t_0)$ to $(x,t)$ and that a simple computation gives $k < \frac{|x_0-x|^2}{t_0-t}$.
\hfill $\square$

\medskip

The set $\big\{ (x_0,t_0), (x_1,t_1), \dots (x_k,t_k), (x,t) \big\}$ appearing in the above proof is often referred to
as \emph{Harnack chain}. By using the following property of the fundamental solution $\Gamma$ of the differential
operator appearing in \eqref{e-unif-par}
\begin{equation}\label{diag}
  \Gamma(0,t) \ge \frac{C}{t^{N/2}}, \qquad \text{for every} \ t>0,
\end{equation}
for some positive constant $C = C(\lambda, \Lambda, N)$. We refer to Nash \cite{Nash} and to Fabes-Strook \cite[Lemma 2.6]{FaSt} for a derivation of \eqref{diag}. By choosing $c = \frac12$ in Theorem
\ref{th-Harnack-NL}, we conclude that there exist two positive constants $C^-, c^-$ such that
\begin{equation*}
  \Gamma(x,t, y,s) \ge \frac{C^-}{(t-s)^{N/2}} \exp\left( - c^- \frac{|x-y|^2}{t-s}\right),
\end{equation*}
for every $(x,t),  (y,s) \in \R^{N+1}$ with $0 < s<t<T$.

\medskip

We explicitly note that the method described above also applies to non-divergence uniformly parabolic operators, if we rely on the Harnack inequality proved by Krylov and Safonov \cite{KrylovSafonov2}. In this setting the inequality \eqref{diag} holds for $t$ belonging to any bounded interval $]0,T[$ and the constant $C$ may depend on $T$. We refer to the manuscript of Konakov \cite{Konakov} for the derivation of \eqref{diag} by using the a parametrix expansion and to the monograph of Bass \cite{Bass} for uniformly parabolic operators with bounded measurable coefficients.

\begin{remark} \label{rem-inv} \ Before considering degenerate parabolic operators, we point out that the method used
in the proof of Theorem \ref{th-Harnack-NL} only relies on the following two ingredients.
\begin{description}
  \item [{\it i)}] The invariance with respect to the Euclidean translation and to the parabolic dilation $(x,t)
\mapsto (x_0 + \rho x, t_0 + \rho^2 t)$ are the properties that allows us to obtain Corollary \ref{cor-Harnack-P} from
Theorem \ref{th-Harnack-P}.
  \item [{\it ii)}] Segments are very simple supports for the construction of Harnack chains. In the study of
degenerate parabolic operators a more sophisticated construction will be needed.
\end{description}
\end{remark}

\section{Degenerate hypoelliptic operators}\label{sec3}
Consider a linear second order differential operator in the form \eqref{e1}
\begin{equation*}
\L  = \sum_{k=1}^{m}X_{k}^{2}+X_{0}-\p_{t}.
\end{equation*}
satisfying the  H\"{o}rmander condition [H]. We introduce a definition based on the vector fields $X_1,...,X_m,Y$.
\begin{definition} \label{def-adm-path}
We say that $\g$ is an $\L$-\emph{admissible path} starting from $z_0 \in \rnn$ if it is an absolutely continuous
solution of the following ODE
\begin{equation*} 
\begin{split}
    & \dot \g(\t) =  \displaystyle{\sum_{k=1}^{m}}
    \o_k(\t) X_k(\g(\t)) + Y(\g(\t)), \\ 
   & \g(0) =  z_0. 
\end{split}
\end{equation*}
 with $\o_1, \dots, \o_{m} \in L^1([0,T])$.

Let $\Omega$ be an open subset of $\rnn$ and $z_0 \in \Omega$. The \emph{attainable set} of $z_0$ in $\Omega$ is
\begin{align*} 
  \A_{z_0} (\Omega) = \big\{z \in \Omega \mid
  &\mbox{there exists an $\L$-admissible path} \ \gamma \ \mbox{such that} \\ 
  & \qquad \qquad   
   \g(0) =  z_0, \ \gamma(T) = z \ \mbox{and} \ \g(\tau) \in \Omega \ \mbox{for} \ 0 \le \tau \le T
\big\}.
\end{align*}
\end{definition}

The following version of the Harnack inequality is based on the definition of attainable set. It has been introduced
in \cite{CMP,CNP} and in its general form in \cite{KogojPolidoro} for operators in the form \eqref{e1}.

\begin{theorem} \label{th-Harnack-KP}
Let $u$ be a non negative solution of $\L u=0$ in some bounded open set $\Omega \subset \rnn$, and let $z_0 \in
\Omega$. Suppose that ${\text{\rm Int}}\left(\overline{\A_{z_0}(\Omega)}\right) \ne \emptyset$. Then, for every compact
set  $\K \subset {\text{\rm Int}}\left(\overline{\A_{z_0}(\Omega)}\right)$ there exists a positive constant
$C_{\K}$, only depending on $\Omega, \K, z_0$ and $\L$, such that
$$
  \sup_{\K}u(z) \le C_{\K}u(z_0).
$$
\end{theorem}

If the operator $\L$ is also invariant with respect to suitable non Euclidean \emph{translations} and \emph{dilations},
then Theorem \ref{th-Harnack-KP} restores an \emph{invariant Harnack inequality} useful for the construction of Harnack
chains.

\smallskip

 \noindent {\sc Hypothesis [G1]}\ {\it There exists a Lie group $\mathbb{G}= \left(\rnn, \circ \right)$ such
that} $X_{1},\dots,X_{m}, Y$ {\it are left invariant on} $\mathbb{G}$, i.e.: {\it given $\xi \in \R^{N+1}$ and denoting
by} $\ell_{\xi}(z)=\xi \circ z,$ {\it the left translation of $z \in \R^{N+1}$ it holds}
\begin{equation*} 
\begin{split}
  & X_i(u(\ell_{\xi}(z))) =(X_iu)(\ell_{\xi}(z)), \qquad i= 1, \dots, m, \\
  & Y(u(\ell_{\xi}(z))) =(Yu)(\ell_{\xi}(z)),
\end{split}
\end{equation*}
{\it for every smooth function $u$.}

As we will see in the next sections, all the examples listed in the Introduction do satisfy the above assumption, that
replaces the usual invariance with respect to the Euclidean translation. For some operators $\L$ considered in this
note,  the vector fields $X_{1},\dots,X_{m}, Y$  are also invariant with respect to a rescaling property
$\left(\d_{\l}\right)_{\l >0}$ of the Lie group $\mathbb{G}$, which replaces the multiplication by a positive scalar in
a vector space.

\smallskip

\noindent {\sc Hypothesis [G2]}\ {\it There exists a dilation $\left(\d_{\l}\right)_{\l >0}$ on the Lie group
$\mathbb{G}$ such that the vector fields $X_{1},\dots,X_{m}$ are $\d_{\l}$-homogeneous of
degree one and $Y$ is $\d_{\l}$-homogeneous of degree two.} i.e.:
\begin{equation*} 
\begin{split}
X_i(u(\d_{\l}(z))) & =\l(X_iu)(\d_{\l}(z)), \qquad i= 1, \dots, m, \\
Y(u(\d_{\l}(z))) & =\l^2(Yu)(\d_{\l}(z)),
\end{split}
\end{equation*}
{\it for every smooth function $u$.}

When both of assumptions [G1] and [G2] are satisfied, we say that $$\mathbb{G} = \left(\rnn,\circ,
\left(\d_{\l}\right)_{\l >0} \right)$$ is a \emph{homogeneous} Lie group and the operator $\L$ is invariant with respect
to the left translations of $\mathbb{G}$, and homogeneous of degree $2$ with respect to the dilation
of $\mathbb{G}$. In this case we easily obtain from Theorem \ref{th-Harnack-KP} an invariant Harnack inequality
analogous to Corollary \ref{cor-Harnack-P}. Consider any bounded open set $\Omega \subset \rnn$ with $0 \in
\Omega$ and suppose that it is \emph{star-shaped} with respect to $\left(\d_{\l}\right)_{\l >0}$, that is
\begin{equation*}
  \delta_r(\Omega) : = \big\{ \delta_r (z) \mid z \in \Omega \big\}\subset \Omega, \qquad \text{for every} \ r
\in ]0,1].
\end{equation*}
If ${\text{\rm Int}}\left(\overline{\A_{0}(\Omega)}\right) \ne \emptyset$, we choose any compact set $\K \subset
{\text{\rm Int}}\left(\overline{\A_{0}(\Omega)}\right)$. For every $r>0$ and $z_0 \in \rnn$ we set
\begin{equation*}
  \Omega_r(z_0)=z_0 \circ \delta_r(\Omega) := \big\{ z_0 \circ \delta_r (z) \mid z \in \Omega \big\}.
\end{equation*}
Note that we also have $z_0 \circ \delta_\rho(\K) \subset {\text{\rm Int}}
\left(\overline{\A_{z_0}(\Omega_r(z_0))}\right)$ for every $\rho \in ]0, r]$, since $\Omega$ is star-shaped. We define
\begin{equation*}
 \P_r(z_0)= \bigcup_{0 < \rho \le r} z_0 \circ \delta_\rho(\K).
\end{equation*}

\begin{theorem} \label{cor-Harnack-KP}
Let $\L$ be an operator in the form \eqref{e1} satisfying assumptions [G1] and [G2] and let $\Omega_r(z_0)$ as above.
Suppose that ${\text{\rm Int}}\left(\overline{\A_{z_0}(\Omega_r(z_0))}\right) \ne \emptyset$, then
\begin{equation*}
  \sup_{{\P_r(z_0)}} u(x,t) \le C_{\K} u (z_0)
\end{equation*}
for every positive solution $u$ of $\L u=0$ in $\Omega_r(z_0)$. Here $C_{\K}$ is the same constant appearing
in Theorem \ref{th-Harnack-KP}.
\end{theorem}

Theorem \ref{cor-Harnack-KP} is the Harnack inequality that replaces Corollary \ref{cor-Harnack-P} in the non Euclidean
setting that is natural for the study of degenerate operators $\L$. In accordance with Remark \ref{rem-inv}, this is
the first ingredient for the construction of Harnack chains. It turns out that the second ingredient is the
$\L$-admissible path, which is the natural substitute of the segment used in the Euclidean setting. To replicate the
construction made in the proof of Theorem \ref{th-Harnack-NL} we only need to choose $\g$, with $\g(0) = (x_0,t_0)$,
and $\P_{(x_0, t_0)}$ with the following property:
\begin{equation} \label{eq-property}
  \text{there exists $s_0 \in ]0,t_0-t[$ such that $\g(s) \in \P_{(x_0, t_0)}$ for $s \in ]0, s_0]$.}
\end{equation}
All the examples in this note satisfy \eqref{eq-property}. Thus we have what we need to construct a Harnack
chain $\big\{ (x_0,t_0), (x_1,t_1), \dots (x_k,t_k), (x,t) \big\}$ with starting point at $(x_0,t_0)$ and end point at
$(x,t)$.

In order to find an accurate bound of the positive solutions of $\L u = 0$ we need to control the \emph{length} $k$ of
the Harnack chain. It is possible to prove that there exists a positive constant $h$ such that, if we construct the
Harnack chain by using the $\L$-admissible path $\g$ as in Definition \ref{def-adm-path}, with $z_0= (x_0,t_0)$ and
$z=(x,t)$, then $T = t_0-t$ and we have
\begin{equation}\label{ecost}
  k \le  \tfrac{1}{h} \Phi(\o) + 1, \qquad \Phi(\o) := \int_0^{t_0-t}\|\o(s)\|^2 d\, s.
\end{equation}
In the sequel we will refer to the integral appearing in \eqref{ecost} as the \emph{cost} of the path $\gamma$
associated to the \emph{control} $(\o_1,...,\o_m)$. We then conclude that there exist three positive constants $\theta,
h$ and $M$, with $\theta < 1$ and $M>1$, only depending on the operator $\L$ such that
\begin{equation}\label{Harnack_chains_result}
 u(x,t) \leq  M^{1+ \frac{\Phi(\o)}{h}} u(x_0, t_0),
\end{equation}
for every positive solution $u$ of $\L u=0$, were $(x_0, t_0), (x,t) \in \rn \times ]0,T[$ are such that $0<t_0-t<\theta
t_0$.

Note that \eqref{Harnack_chains_result} provides us with a bound depending on the choice of the $\L$-admissible path
$\gamma$ steering $(x_0, t_0)$ to $(x,t)$. In order to get the best exponent, we can optimize the choice of $\g$.  With
this spirit, we define the  \emph{Value function}
\begin{equation}\label{Value_f}
 \Psi(x_0, t_0;x,t) = \inf_{\o} \big\{\Phi(\o) \big\},
\end{equation}
where the \emph{infimum} is taken in the set of all the $\L$-admissible paths $\gamma$ steering $(x_0, t_0)$ to
$(x,t)$, and satisfying \eqref{eq-property}. We summarize this construction in the following general statement.

\medskip

\noindent {\it
Let $\L$ be an operator in the form \eqref{e1} satisfying conditions [H], [G1] and [G2], and assume that there is a
positive $r$ and an open star-shaped set $\Omega$ with $0 \in \Omega$ such that ${\text{\rm Int}}\left(
\overline{\A_{0}(\Omega_r(0))}\right) \ne \emptyset$.
Moreover, if all the $\L$-admissible paths $\gamma$ steering
$(x_0, t_0)$ to $(x,t)$ satisfy \eqref{eq-property}, then there exist three positive constants $\theta,
h$ and $M$, with $\theta < 1$ and $M>1$, only depending on the operator $\L$ such that the following property holds.}

{\it Let $(x_0, t_0), (x,t) \in \rnn$ with $0<t_0-t<\theta t_0$. Then, for every positive solution $u: \rn \times
]0,T[$ of $\L u=0$ it holds}
\begin{equation} \label{Harnack_chains}
 u(x,t) \leq  M^{1+ \frac{1}{h}\Psi(x_0, t_0;x,t)} u(x_0, t_0).
\end{equation}

Inequality \eqref{Harnack_chains}
is the main step in the proof of our lower bound for the fundamental solution. All the
examples considered in this note satisfy conditions [H], [G1]. Some examples also satisfy [G2], some examples do not.
However, in this case, a scale invariant Harnack inequality still holds true, then the method provides us with a lower
bound of the fundamental solution.

\section{Degenerate hypoelliptic operators on homogeneous groups}\label{sec4}
The Heat operator on the Heisenberg group
$$\L =  X_1^2 +X_2^2  - \p_t
$$
where
$$
  X_1=\p_x - \tfrac{1}{2}y \p_w, \qquad X_2=\p_y + \tfrac{1}{2}x \p_w
$$
are vector fields acting on the variable $(x,y,w,t) \in \R^{4}$, is the simplest example of degenerate operator built
by a sub-Laplacian on a stratified Lie group. The vector fields $X_1, X_2$ are invariant with respect to the left
translation on the Heisenberg group on $\R^3$, whose operation is defined as
\begin{equation*}
  (x_0,y_0,w_0) \circ (x,y,w) = \left(x_0 + x,y_0+y,w_0+w + \tfrac12 (x_0 y - y_0 x)\right).
\end{equation*}
The above operation is extended to $\R^4$ by setting
\begin{equation*}
  (x_0,y_0,w_0,t_0) \circ (x,y,w,t) = \left(x_0 + x,y_0+y,w_0+w + \tfrac12 (x_0 y - y_0 x), t_0+ t\right).
\end{equation*}
Moreover $\L$ is invariant with respect to the following dilation
\begin{equation*}
  \delta_r (x,y,w,t) = \left(r x,r y ,r^2 w , r^2 t \right),
\end{equation*}
then the hypotheses [G1] and [G2] are fulfilled by $\L$. Furthermore, it satisfies the following property.

\smallskip

\noindent [C] {\it For every $x_0, x \in\rn$, and for every positive $T$ there exists an absolutely continuous path
$\g_0: [0,T] \to \rn$ such that}
\begin{equation}\label{eq-dist}
  \dot \g_0 (\t) =  \displaystyle{\sum_{k=1}^{m}} \o_k(\t) X_k(\g_0(\t)), \quad  \g_0 (0) =  x_0, \ \g_0(T) =  x.
\end{equation}

\smallskip

Note that, for operators $\L$ in the form \eqref{e1} with $X_0=0$, condition [C] is equivalent to the \emph{strong} H\"{o}rmander condition
\begin{equation*}
 \text{\rm rank Lie}\big\{X_{1},\dots,X_{m}\big\}(x)=N, \qquad
 \forall x\in \R^{N}.
\end{equation*}
Moreover, for every $\Omega \subset \rnn$ and for every $(x_0,t_0) \in \Omega$ there exist a positive $\varepsilon$ and
 a neighborhood $U$ of $x_0$ such that $U \times ]t_0, t_0- \varepsilon[ \subset  \A_{(x_0,t_0)}(\Omega)$. This
particular geometric property of the attainable set implies that an invariant Harnack inequality analogous to the
standard parabolic one holds for this operator. The only difference is that the Euclidean translation and the parabolic
dilations are replaced by the operations used to satisfy hypotheses [G1] and [G2]. In conclusion, the hypotheses we
need to prove \eqref{Harnack_chains} are satisfied by the heat operator on the Heisenberg group. In particular, this
method leads us to the lower bound of the following version of \eqref{e-keystone}: there exist positive constants $c^-,
C^-, c^+,C^+$ such that
\begin{equation} \label{e-keystone-CC}
  \tfrac{c^-}{\sqrt{
  |\B_{t-\tau}(x)|}} \exp \left(- C^- \tfrac{d_{CC}(x,\xi)^2}{t-\tau} \right) \le \Gamma(x,t,\xi,\tau) \le
  \tfrac{C^+}{\sqrt{
  |\B_{t-\tau}(x)|}} \exp \left(- c^+ \tfrac{d_{CC}(x,\xi)^2}{t-\tau} \right),
\end{equation}
where $d_{CC}$ denotes the \emph{Carnot-Caratheodory distance}
$$
  d_{CC}(x_0, x) = \inf \{ \ell(\g_0)\  | \ \g_0 \ \text{is as in \eqref{eq-dist}}\}, \quad
  \ell(\g):=\int_0^{T}\|\o(s)\|ds.
$$
and $|\B_{r}(x)|$ is the volume of the metric ball with center at $x$ and
radius $r$. To make more precise the analogy between \eqref{e-keystone} and \eqref{e-keystone-CC}, we recall that if $\mathbb{H}$ is a homogeneous Lie group on $\rn$, then
$$|\B_{r}(x)|=r^Q|\B_{1}(0)|,$$
where $Q$ is an integer called \emph{homogeneous dimension} of $\mathbb{H}$.
We recall that the upper bound was proved by Davies in \cite{Davies}, and the upper and lower bounds are due to Jerison
and S{\'a}nchez-Calle  \cite{Jerison} and to Varopoulos, Saloff-Coste and Coulhon \cite{VSC}. Note that
$\Psi(x_0,t_0; x,t) = \frac{d_{CC}(x_0, x)^2}{t_0 - t}$. Indeed,  if we consider the path $\gamma (s) = (\gamma_0(s),
t_0 -s)$ with $0 \le s \le t_0 -t$, then by the Cauchy-Schwartz inequality, we obtain $\ell(\g_0) \le \sqrt{\Phi(\omega)}
\sqrt{ t_0 -t}$. Moreover the equality occurs only if the norm of the control $\o$ is constant, that is
\begin{equation*}
  \ell(\g_0) = \sqrt{\Phi(\omega)} \sqrt{ t_0 -t} \iff (\o_1^2+...+\o_m^2)(s) = \frac{\Phi(\o)}{ t_0 -t} \ \text{for
every} \ s \in [0,  t_0 -t].
\end{equation*}
We refer to the article \cite{BoscainPolidoro} for the study of a more general class of operator satisfying [G1], [G2]
and [C], that includes heat operators on Carnot groups and also operators $\L$ with $X_0 \ne 0$. We also recall that in
the article \cite{CintiPolidoro} the analogous upper bound has been proved by using a PDE method combined with the
Optimal Control Theory.

\section{Degenerate Kolmogorov equations}\label{sec5}
The simplest degenerate example of degenerate Kolmogorov operator is
\begin{equation}\label{Kolmogorov_op}
\L :=  \p_{x}^2 + x \p_y - \partial_t, \qquad (x,y,t) \in \R^2\times ]0,T[,
\end{equation}
it writes in the form \eqref{e1}, if the vector fields $X,\, Y$ are
      $$
X(x,y,t) = \partial_x  \sim \left(
                  \begin{array}{c}
                    1 \\
                    0 \\
                    0 \\
                  \end{array}
                \right),
                \qquad Y(x,y,t) = x \partial_y - \partial_t \sim  \left(
                  \begin{array}{c}
                    0 \\
                    x \\
                    -1 \\
                  \end{array}
                \right).
 $$
$\L$ is related to the following  stochastic process
\begin{equation} \label{Langevin}
\begin{cases}
  & X_t = x_0 + W_t, \\
  & Y_t = y_0 + \int_0^t (x_0 + W_s) \, d s.
\end{cases}
\end{equation}
which satisfies the \emph{Langevin equation} $dX_t = dW_t, dY_t = X_t dt$. We recall that this kind of stochastic
process appears in several research areas. For instance, in Kinetic Theory, $(X_t)_{t \ge 0}$ describes the velocity of
a particle, while $(Y_t)_{t \ge 0}$ is its position. We note that
\begin{description}
\item [{\it i)}] $X$ and $Y$ are invariant with respect to the left translation of the group defined by the following operation  
\begin{equation} \label{transl_Kolm}
(x_0,y_0,t_0) \circ (x,y,t)= (x+x_0, y + y_0 - tx_0, t+t_0), \qquad (x,y,t), (x_0,y_0,t_0) \in \mathbb{R}^3,
\end{equation}
\item [{\it ii)}] $X$ and $Y$ are  homogeneous of degree $1$ and $2$, respectively, with respect to the dilation
\begin{equation} \label{dil_Kolm}
(\delta_{\rho})_{\rho > 0}: (x,y,t) \mapsto (\rho x, \rho^3 y, \rho^2 t)= {\rm diag}(\rho , \rho^3 , \rho^2) \cdot \left(
                 \begin{array}{c}
                   x \\
                   y \\
                   t \\
                 \end{array}
               \right).
\end{equation}
In particular, $\L$ satisfies the Hypotheses [G1] and [G2].
  \item [{\it iii)}] The $\L$-admissible paths are the solutions $\g(s)=(x(s),y(s),t(s))$ of  the following equation
\begin{equation*} 
 \left\{
   \begin{array}{ll}
     \dot{x}(s)=\o(s), & \hbox{$x(0)=x_0$,} \\
     \dot{y}(s)=x(s), & \hbox{$y(0)=y_0$,} \\
     \dot{t}(s)=-1, & \hbox{$t(0)=t_0$.}
   \end{array}
 \right.
\end{equation*}
\end{description}
It is easy to check that the attainable set of the point $(0,0,0)$ in the open set $\Omega = ]-1,1[^3$ is
$\mathcal{A}_{(0,0,0)}(\Omega) = \big\{ (x,y,t) \in \Omega \mid t< - | y| \big\},$ (see Fig. 4).
\begin{center}
 \begin{pspicture*}(-5,-2.5)(5,1.5)
\scalebox{.80}{%
\pcline[linecolor=dgrey](0,-2)(0,1.4)
\pcline[linecolor=dgrey](0,1.4)(-.05,1.2)%
\pcline[linecolor=dgrey](0,1.4)(.05,1.2)%
\uput[0](0,1.4){$t$}
\pcline[linecolor=dgrey](-4,1)(4,-1)%
\pcline[linecolor=dgrey](4,-1)(3.8,-1)%
\pcline[linecolor=dgrey](4,-1)(3.9,-.9)%
\uput[0](4,-1){$x$}
\pcline[linecolor=dgrey](2,1)(-2,-1)%
\pcline[linecolor=dgrey](-2,-1)(-1.8,-1)%
\pcline[linecolor=dgrey](-2,-1)(-1.9,-.9)%
\uput[180](-2,-1){$y$}
\dotnode(0,0){O}\nput{115}{O}{$(0,0,0)$}
\pcline[linecolor=dgrey](-3,0)(1,-1)%
\pcline[linecolor=dgrey](-3,-2)(1,-3)%
\pcline[linecolor=dgrey](3,0)(-1,1)%
\pcline[linecolor=dgrey, linestyle=dashed](3,-2)(-1,-1)%
\pcline[linecolor=dgrey](3,0)(3,-2)%
\pcline[linecolor=dgrey](-3,0)(-3,-2)%
\pcline[linecolor=dgrey](1,-1)(1,-3)%
\pcline[linecolor=dgrey, linestyle=dashed](-1,1)(-1,-1)%
\pcline[linecolor=dgrey](3,0)(1,-1)%
\pcline[linecolor=dgrey](-3,0)(-1,1)%
\pcline[linecolor=dgrey](3,-2)(1,-3)%
\pcline[linecolor=dgrey, linestyle=dashed](-3,-2)(-1,-1)%
\pcline[linecolor=dgrey](-2,.5)(2,-.5)%
\pcline[linecolor=dgrey, linestyle=dashed](0,0)(0,-2)%
\pcline[linecolor=dgrey](2,-.5)(1,-3)%
\pcline[linecolor=dgrey, linestyle=dashed](-2,.5)(-1,-1)%
\pspolygon[fillstyle=solid,fillcolor=lgrey](-2,.5)(2,-.5)(1,-3)(-3,-2)
\pspolygon[fillstyle=solid,fillcolor=grey](2,-.5)(3,-2)(1,-3)
\pcline[linecolor=dgrey](2,1)(-2,-1)%
\pcline[linecolor=dgrey](-2,-1)(-1.8,-1)%
\pcline[linecolor=dgrey](-2,-1)(-1.9,-.9)%
\uput[180](-2,-1){$y$}
\dotnode(0,0){O}\nput{115}{O}{$(0,0,0)$}
\pcline[linecolor=dgrey](-3,0)(1,-1)%
\pcline[linecolor=dgrey](-3,-2)(1,-3)%
\pcline[linecolor=dgrey](3,0)(-1,1)%
\pcline[linecolor=dgrey, linestyle=dashed](3,-2)(-1,-1)%
\pcline[linecolor=dgrey](3,0)(3,-2)%
\pcline[linecolor=dgrey](-3,0)(-3,-2)%
\pcline[linecolor=dgrey](1,-1)(1,-3)%
\pcline[linecolor=dgrey, linestyle=dashed](-1,1)(-1,-1)%
\pcline[linecolor=dgrey](3,0)(1,-1)%
\pcline[linecolor=dgrey](-3,0)(-1,1)%
\pcline[linecolor=dgrey](3,-2)(1,-3)%
\pcline[linecolor=dgrey, linestyle=dashed](-3,-2)(-1,-1)%
\pcline[linecolor=dgrey](-2,.5)(2,-.5)%
\pcline[linecolor=dgrey, linestyle=dashed](0,0)(0,-2)%
\pcline[linecolor=dgrey](2,-.5)(1,-3)%
\pcline[linecolor=dgrey, linestyle=dashed](-2,.5)(-1,-1)%
}
\end{pspicture*}
\end{center}
\begin{center}
  {\scriptsize \sc Fig. 4 - $\mathcal{A}_{(0,0,0)}(\Omega)$.}
\end{center}

As the interior of $\mathcal{A}_{(0,0,0)}(\Omega)$ is not empty, Theorem \ref{cor-Harnack-KP} gives an invariant
Harnack inequality for $\L$, and we can apply \eqref{Harnack_chains} to prove lower bounds for positive solutions
defined on the domain $\R^2 \times ]0,T[$. The Optimal Control Theory provides us with an explicit expression of the
value function $\Psi_0$ for $\L$ in \eqref{Kolmogorov_op}
\begin{equation} \label{Value_kolm}
  \Psi_0( x,y,t;\xi, \eta, \tau) = \frac{(x-\xi)^2}{t- \tau} + \frac{12}{(t-\tau)^3} \left(y - \eta - (t-\tau)
\tfrac{(x+\xi)}{2}\right)^2.
\end{equation}
This is a remarkable fact, as it is known that the explicit expression of the fundamental solution of $\L$ was written
by Kolmogorov (1934) and is
\begin{equation} \label{Fund_sol}
 \Gamma_0 (x,y,t; \xi, \eta,\tau) = \frac{ \sqrt{3}}{ 2\pi (t-\tau)^2} \exp \left( - \tfrac{(x-\xi)^2}{4(t- \tau)} -
\tfrac{3}{(t-\tau)^3} \left(y - \eta - (t-\tau) \tfrac{(x+\xi)}{2}\right)^2 \right).
\end{equation}

We briefly discuss here the anisotropic dilation \eqref{dil_Kolm}. We first note that the H\"{o}rmander condition is satisfied since
$$[X,Y]=XY-YX =\p_y \sim \left(
                           \begin{array}{c}
                             0 \\
                             1 \\
                             0 \\
                           \end{array}
                         \right)
$$
and that $\p_y$ is homogeneous of degree three as $XY$ and $YX$ are both homogeneous of degree three. This explain the exponent $3$ appearing in \eqref{dil_Kolm}. Moreover, since $${\rm det}\big( {\rm diag}(\rho,\rho^3)\big)=\rho^4,$$ then $Q=4$ is the spatial homogeneous dimension of $\R^2$ with respect to the dilation \eqref{dil_Kolm}. Furthermore, in view of \eqref{Langevin}, such dilation has a natural probabilistic meaning as one has $\text{Var}(X_t)=t$ and $\text{Var}(Y_t)=t^3/3$.

\medskip

The lower bound based on the value function $\Psi$ is useful as we consider Kolmogorov equations in the form
\begin{equation} \label{e-Kolmo-gen}
  \partial_{t} u (x,t) = \sum_{i,j=1}^m a_{ij}(x,t) \partial^2_{x_ix_j} u(x,t) +
\sum_{i,j=1}^N b_{ij}x_j \partial_{x_i}u(x,t), \qquad (x,t)
\in \R^N \times ]0,T[,
\end{equation}
with bounded H\"older continuous coefficients $a_{ij}$'s. In the study of this family of operators, we assume that $m
<N$, the matrix $\left(a_{ij}(t,x)\right)_{i,j = 1 ,\dots, m}$ is uniformly positive in $\R^m$. Moreover, the
H\"ormander condition is satisfied for the operator $\L_{(\xi,\tau)}$ \emph{frozen} at some point $(\xi,\tau) \in \rnn$,
that is obtained from the equation in  \eqref{e-Kolmo-gen} by replacing every function $a_{ij}= a_{ij}(x,t)$ with
$a_{ij}(\xi,\tau)$. It turns out that this condition does not depend on the choice of the point $(\xi,\tau)$, that
$\L_{(\xi,\tau)}$ is invariant with respect to a Lie group $\mathbb {G}$ on $\rnn$ which does not depend on $(\xi,\tau)$.
In this case the \emph{parametrix method} provides us with the existence of a fundamental solution $\Gamma$ of the
operator introduced in \eqref{e-Kolmo-gen}. The method also gives an upper bound of the form
\begin{equation*} 
 \Gamma (x,t;\xi,\tau) \le \frac{C^+}{(t-\tau)^{Q/2}} \exp \left( - c^+ \Psi(x,t;\xi,\tau) \right) \quad (\xi,\tau), (x,t)
\in \R^N \times ]0,T[, \ \ t>\tau,
\end{equation*}
where $Q$ is the \emph{homogeneous dimension} of the space $\R^N$ with respect to the underlying Lie Group in $\rnn$,
and  $C^+, c^+$ are constants depending on the operator.  The method described in this section gives the
analogous lower bound for $\Gamma$
\begin{equation*} 
 \frac{c^-}{(t-t_0)^{Q/2}} \exp \left( - C^- \Psi(x,t;x_0,t_0) \right) \le \Gamma (x,t;x_0,t_0)  \qquad (x_0,t_0),
(x,t) \in \R^N \times ]0,T[.
\end{equation*}


We conclude this section with a discussion on another meaningful example of operator which writes in the form
\eqref{e-Kolmo-gen} and is somehow more degenerate than \eqref{Kolmogorov_op}. It is
\begin{equation}\label{Kolm_general}
\L=\p_{x_1}^2+x_1 \p_{x_2}+....+x_{N-1}\p_{x_N}-\p_t,
\end{equation}
which is related to the following stochastic process
\begin{equation}\label{iterated}
dX_t^1=dW_t, \quad dX_t^2=X_t^1dt, \quad .... \quad, dX_t^N=X^{N-1}_t dt, \quad t \ge 0.
\end{equation}
As the operator defined in \eqref{Kolmogorov_op}, the one in \eqref{Kolm_general} can be written as $\L=X^2+Y$ with:
$$
X(x,t) = \partial_{x_1}  \sim \left(
                  \begin{array}{c}
                    1 \\
                    0 \\

                   \vdots\\
                    0 \\
                  \end{array}
                \right),
                \qquad Y(x,t) = \sum_{j=1}^{N-1}x_j \partial_{x_{j+1}} - \partial_t \sim  \left(
                  \begin{array}{c}
                    0 \\
                    x_1 \\
                    x_2\\
                    \vdots\\
                    -1 \\
                  \end{array}
                \right).
 $$
Note that, in this case, $\p_{x_{j+1}}=[\p_{x_j},Y]$ for $j=1,...,N-1$.
As a consequence, $\L$ is invariant with respect to the dilation defined by the following matrix:
$${\rm diag}(\rho ,\rho^3,...,\rho^{2N-1},\rho^2),$$ then its homogeneous dimension $Q$ is equal to
$N^2$. Accordingly, we have that $\text{Var}(X_t^j)=c_j t^{2j-1}, \ j=1,...,N$, where $c_j$ is a positive constant.

\medskip

We recall that the parametrix method has been used by several authors for the study of degenerate Kolmogorov equations.
We recall the works of Weber \cite{Weber}, Il'In \cite{Il'in}, Sonin \cite{Sonin}, Polidoro \cite{Polidoro2, Polidoro1},
Di Francesco and Polidoro \cite{DiFrancescoPolidoro}. In particular, the lower bound of the fundamental is proved in
\cite{Polidoro1} and in \cite{DiFrancescoPolidoro}.

More recently, Delarue and Menozzi \cite{DelarueMenozzi} extended the above bounds to a class of Degenerate Kolmogorov Operator with possibly non-linear drifts satisfying H\"{o}rmander condition, under spatial H\"{o}lder continuity  assumptions on the coefficients
$a_{ij}$'s. They obtained analogous bounds by combining stochastic control methods with the parametrix representation of
the fundamental solution given by McKean and Singer in \cite{Singer}.

\section{More degenerate equations}\label{sec6}
In this section we consider a stochastic process studied By Cinti, Menozzi and Polidoro in
\cite{CMP}. It is similar to the one considered in Section \ref{sec4}, as it writes as follows
\begin{equation} \label{eq-CMP}
  \L :=  \p_{x}^2 + x^2 \p_y - \partial_t, \qquad (x,y,t) \in \R^2\times (0,T),
\end{equation}
and is related to the following stochastic differential equation
\begin{equation} \label{eq-S-CMP}
\begin{cases}
  & X_t = x_0 + W_t, \\
  & Y_t = y_0 + \int_0^t (x_0 + W_s)^2 \, d s.
\end{cases}
\end{equation}
A representation of the density of this process has been obtained from the seminal works of Kac \cite{Kac} in terms of
the Laplace transform of the process $(Y_t)_{t \ge 0}$. We also refer to the monograph of Borodin and Salminen
\cite{boro:salm:02} for an expression in terms of special functions. We also quote the works of  Smirnov \cite{Smirnov}
and Tolmatz \cite{Tolmatz} on the distribution function of the square of the Brownian bridge.

We give explicit upper and lower bounds for the density of the process $(X_t,Y_t)_{t \ge 0}$ by the approach described
in Section \ref{sec3}. Note that new difficulties appear in the study of the operator $\L$ defined in \eqref{eq-CMP}.
Indeed,
if we  write $\L$ as follows
\begin{equation*}
  \L  = X^2 + Y, \quad  \text{with} \quad X = \partial_x, \ Y = x^2 \p_y - \partial_t,
\end{equation*}
then the commutator $[X,Y] (x,y,t) = 2 x \partial_y$ vanishes in the set $\big\{ x=0 \big\}$, and we need a second
commutator $[X, [X,Y]] (x,y,t) = 2 \partial_y$ to satisfy the H\"ormander condition at every point of $\R^3$. As a
consequence, a Lie group leaving invariant the equation $\L u = 0$ cannot exist. This problem is overcome by a
\emph{lifting procedure} (see Rothshild and Stein \cite{RothschildStein}). Specifically, we consider the following operator
\begin{equation*} 
  \widetilde {\L} :=  \p_{x}^2 + x \partial_w + x^2 \p_y - \partial_t, \qquad (x,y,w,t) \in \R^3\times (0,T),
\end{equation*}
and we consider any solution of $\L u = 0$ as a function that does not depend on $w$, and that solves the equation
$\widetilde {\L} u = 0$. The lifting procedure allows us to rely on the Lie group invariance of $\widetilde {\L}$ in
the study of the positive solutions of $\L u = 0$. Indeed, we have
\begin{description}
\item [{\it i)}] The operator $\widetilde {\L}$ is invariant with respect to the following Lie group operation
\begin{equation*} 
  (x_0,y_0,w_0, t_0) \circ (x,y,w,t)= (x+x_0, y + y_0 + 2 x_0 w - t x_0^2, w + w_0 - t x_0, t+t_0),
\end{equation*}
defined for every  $(x,y,w,t), (x_0,y_0, w_0, t_0) \in \mathbb{R}^4$. In particular, it holds
\begin{equation*}
  (\widetilde \L u) (z_0 \circ z) = \widetilde \L(u(z_0\circ z)),
\end{equation*}
for every $z_0=(x_0,y_0,w_0,t_0)$ and $z=(x,y,w,t) \in \R^4$.
\item [{\it ii)}] The operator $\widetilde \L$ is invariant with respect to the following dilation
 \begin{equation*} 
  (\delta_{\rho})_{\rho \ge 0}: (x,y,w,t) \mapsto (\rho x, \rho^4 y, \rho^3 w, \rho^2 t).
 \end{equation*}
That is, it holds:
\begin{equation*}
  \rho^2 \, (\L u)(\rho x, \rho^3 y, \rho^2 t)=\L(u(\rho x, \rho^3 y, \rho^2 t)).
\end{equation*}
  \item [{\it iii)}] The attainable set of the \emph{origin} in the box $\Omega = ]-1,1[^4$ is
\begin{equation*}
\mathcal{A}_{(0,0,0,0)} (\Omega) =\Big\{ (x,w,y,t) \in ]-1,1[^4 \ \mid 0 \le y \le -t, w^2 \le - t y \Big\}.
\end{equation*}
Figure 4 describes the projection on the hyperplane $\big\{ x= 0 \big\}$ of the set $\mathcal{A}_{(0,0,0,0)}$
\begin{center}
\begin{pspicture*}(-3.1,-2.5)(5,1.5)
\scalebox{.80}{%
\begin{psclip}{\pspolygon[linecolor=dgrey](-3,0)(-1,1)(3,0)(3,-2)(1,-3)(-3,-2)(-3,0)}%
\pspolygon[linecolor=dgrey, fillstyle=solid,fillcolor=grey]
(2,-3.2)(1,-3)(.5,-2.82)(0,-2.55)(-.26,-2.3)(-.36,-2.1)(-.34,-1.9)(0,0)(3,-2)
\pspolygon[linecolor=dgrey, fillstyle=solid, fillcolor=lgrey] (1,-3)(0,0)(3,-2)
\psbezier[linecolor=dgrey, linestyle=dashed](-.33,-2.1)(-.3,-1.2)(2.2,-2)(3,-2)
\end{psclip}
\pcline[linecolor=dgrey](0,-2)(0,1.4)
\pcline[linecolor=dgrey](0,1.4)(-.05,1.2)%
\pcline[linecolor=dgrey](0,1.4)(.05,1.2)%
\uput[0](0,1.4){$t$}
\pcline[linecolor=dgrey](-4,1)(4,-1)%
\pcline[linecolor=dgrey](4,-1)(3.8,-1)%
\pcline[linecolor=dgrey](4,-1)(3.9,-.9)%
\uput[0](4,-1){$y$}
\pcline[linecolor=dgrey](-3,0)(1,-1)%
\pcline[linecolor=dgrey](-3,-2)(1,-3)%
\pcline[linecolor=dgrey](3,0)(-1,1)%
\pcline[linecolor=dgrey, linestyle=dashed](3,-2)(-1,-1)%
\pcline[linecolor=dgrey](3,0)(3,-2)%
\pcline[linecolor=dgrey](-3,0)(-3,-2)%
\pcline[linecolor=dgrey](1,-1)(1,-3)%
\pcline[linecolor=dgrey, linestyle=dashed](-1,1)(-1,-1)%
\pcline[linecolor=dgrey](3,0)(1,-1)%
\pcline[linecolor=dgrey](-3,0)(-1,1)%
\pcline[linecolor=dgrey](3,-2)(1,-3)%
\pcline[linecolor=dgrey, linestyle=dashed](-3,-2)(-1,-1)%
\pcline[linecolor=dgrey](2,1)(-2,-1)%
\pcline[linecolor=dgrey](-2,-1)(-1.8,-1)%
\pcline[linecolor=dgrey](-2,-1)(-1.9,-.9)%
\uput[180](-2,-1){$w$}
\dotnode(0,0){O}\nput{115}{O}{$(0,0,0)$}
\uput[180] (4,1) {$\big\{ x= 0 \big\}$}
}
\end{pspicture*}
\end{center}
\begin{center}
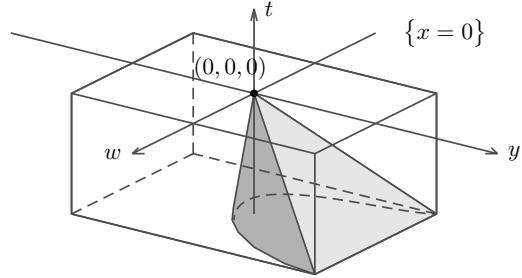
  {\scriptsize \sc Fig. 5 - Projection  of $\mathcal{A}_{(0,0,0,0)} (\Omega)$ on the set $\big\{ x= 0 \big\}$.}
\end{center}
\end{description}

Then, an invariant Harnack inequality needed to construct Harnack chains for the positive solutions of $\widetilde{\L} u
= 0$ is available. The main result of the article \cite{CMP} is the following

\begin{theorem} \label{th-CMP}
Let $\Gamma$ denote the fundamental solution of $\p_{xx} + x^2 \p_y - \p_t$.
\begin{itemize}
\item If $\eta - y \le 0$, then $\Gamma(x,y,t,\xi, \eta, \tau) = 0$;
\item if $\frac{\eta - y}{(t-\tau)^2} > \frac{x^2+\xi^2}{t - \tau} + 1$, then
$$
  \Gamma(x,y,t,\xi, \eta, \tau) \approx \frac{1}{(t-\tau)^{5/2}} \exp \left( - C \left(\frac{(x-\xi)^2}{t-\tau}  +
\frac{\eta - y }{(t-\tau)^2}\right)\right);
$$
\item if $0 < \frac{\eta - y}{(t-\tau)^2} < \frac{1}{2}$, then
$$
  \Gamma(x,y,t,\xi, \eta, \tau) \approx \frac{1}{(t-\tau)^{5/2}} \exp \left( - C \left(\frac{x^4 +\xi^4 +
(t-\tau)^2}{\eta - y }\right)\right).
$$
\end{itemize}
\end{theorem}

We conclude this section with some remarks. We first note that, because of the particular form of the attainable set
$\mathcal{A}_{(0,0,0,0)} (\Omega)$, it is not true that all the $\L$-admissible paths $\gamma$ steering $z_0$ to $z$
satisfy \eqref{eq-property}. For this reason, in the proof of our main result we do not solve any optimal control
problem. We prove our lower bound by choosing smart admissible paths. This construction does not guarantee the
optimality of the lower bounds. However, the comparison with the upper bound, that has the same asymptotic behavior,
shows the optimality of both of them. The diagonal bounds and the upper bounds have been obtained by using probabilistic
methods, and Malliavin Calculus in particular.

We eventually recall that more general operators and stochastic processes are studied in \cite{CMP}. Precisely, we
consider for every positive integer $k$ the process $(X_t, Y_t)_{t \ge 0}$, with value in $\R^n \times \R$
$$
\begin{cases}
  & X_t = x + W_t \\
  & Y_t = y + \int_0^t \sum_j (x + W_s)_j^{k} \, d s,
\end{cases}
$$
whose Kolmogorov equation is
$$
\L := \tfrac12 \Delta_x + (x_1^{k} + \dots + x_n^{k}) \p_y - \partial_t
$$
and
$$
\begin{cases}
  & X_t = x + W_t, \quad (k \ \text{even}) \\
  & Y_t = y + \int_0^t |x + W_s|^{k} \, d s,
\end{cases}
$$
whose Kolmogorov equation is
$$
\L := \tfrac12 \Delta_x + |x|^{k} \p_y - \partial_t.
$$
We refer to the article \cite{CMP} for the precise statement of our achievements and for further details.

\section{Operators related to Arithmetic Average Asian Options}\label{sec7}
In this section we consider the operator
$$\L=x^2 \p_{xx} +x\p_x+x\p_y-\p_t$$
 with $(x,y,t) \in \R^+\times \R \times (0,T)$. It appears in the Black and Scholes setting when we consider the
pricing problem for Arithmetic Average Asian Option. Specifically, we assume that the price of an asset
$(X_t)_{t \ge 0}$ is described by a Geometric Brownian Motion and that the option depends on the arithmetic average of
$(X_t)_{t \ge 0}$. Then, according to the Black and Scholes theory,  the value of the option $v$ is modeled by a
function $v= v(t, X_t, Y_t)$ where
\begin{equation}\label{Yor process}
\begin{cases}
  &X_t = x_0 e^{\sqrt 2 W_t}, \\
  &Y_t = y_0 + x_0\int_0^t e^{\sqrt 2 W_s}ds.
\end{cases}
\end{equation}
This system was widely studied by Yor who wrote in \cite[Chapter 6]{Yor5} its joint density (see equation (6.e) therein)
\begin{equation}\label{e-Yor-density}
 p(x, y, t;x_0,y_0) = \frac{\sqrt{x_0}}{2 \sqrt x (y-y_0)^2}\frac{e^{\frac{\pi^2}{t}}}{\pi\sqrt{ \pi t}}\exp\left( -\frac{x+x_0}{2(y-y_0)}\right)
 \psi\left(\frac{\sqrt {xx_0}}{y-y_0},\frac{t}{2} \right),
\end{equation}
where
\begin{equation}\label{psi}
  \psi\left(z,t \right)=\int_0^{\infty}e^{-\frac{\xi^2}{2 t}}e^{-z\cosh(\xi)}
  \sinh{(\xi)}\sin\left( \frac{\pi \xi}{t} \right)d \xi.
\end{equation}
As in the previous example, the density of the stochastic process $(X_t,Y_t)_{t \ge 0}$  is not strictly positive in
the whole set $\R^+ \times \R \times (0,T)$. In particular, its support is $\R^+ \times (y_0,+\infty) \times (t_0,T)$.

Monti and Pascucci observe in \cite{PascucciMonti2009} that $\L$ is invariant with respect to the following group
operation on $\R^+ \times \R^2$:
\begin{equation}\label{sx}
  (x_0,y_0,t_0)\circ(x,y,t)=(x_0 x,y_0+x_0 y,t_0+t).
\end{equation}
Indeed, if we set
\begin{equation}\label{eq-u-vs-v}
  v(x,y,t) = u (x_0x,y_0+x_0y,t_0+t),
\end{equation}
then $\L v = 0$ if, and only if $\L u = 0$.

Note that $\L$ is not invariant with respect to any dilation group $(\delta_\rho)_{\rho \ge 0}$. On the other hand, as
\begin{equation*}
  \L =  X^2 + Y, \quad \text{with} \quad X(x,y,t) = x \partial_x, \  Y(x,y,t) = x \partial_y   -\partial_t,
\end{equation*}
we have that $\L$ can be approximated by the Kolmogorov operator \eqref{Kolmogorov_op} defined in Section
\ref{sec5}. Indeed, we can consider the coefficient $x$ of the vector field $X$ as a smooth function that is bounded
and bounded by below on every compact set $K \subset \R^+\times \R \times (0,T)$. For this reason, the Harnack
inequality introduced in Serction \ref{sec5} also applies to $\L$.

The $\L$ admissible paths are the solutions of the following differential equation
\begin{equation*} 
 \left\{
   \begin{array}{ll}
     \dot{x}(s)=\o(s)x(s), & \hbox{$x(0)=x_0$,} \\
     \dot{y}(s)=x(s), & \hbox{$y(0)=y_0$,} \\
     \dot{t}(s)=-1, & \hbox{$t(0)=t_0$,}
   \end{array}
 \right.
  \end{equation*}
and we denote by  $\Psi(x_0,y_0,t_0, x,y,t)$ the value function of the relevant optimal control problem with quadratic
cost. The main result for the fundamental solution $\G(x,y,t;x_0,y_0,t_0)$ of the operator $\L$ is  the following
\begin{theorem} \label{th-main}
Let $\Gamma$ be the fundamental solution of $\L$. Then, for every $(x_0,y_0,t_0) \in \R^+ \times \R \times
[0,T[$ we have
\begin{equation} \label{supp_gamma}
  \Gamma(x,y,t, x_0,y_0,t_0)=0 \qquad \forall \ (x,y,t) \in \R^+\times\R^2 \setminus
  \big\{ ]-\infty,y_0[ \times]t_0,T[ \big\}.
\end{equation}
Moreover, for arbitrary $\eps \in ]0, 1[$, there exist two positive constants $c_{\eps}^-,
C_{\eps}^+$ depending on $\eps$, on $T$ and on the operator $\L$, and two positive constants $C^-,
c^+$, only depending on the operator $\L$ such that
\begin{equation} \label{e-twosidedbounds1}
\begin{split}
  \frac{ c_{\eps}^-}{x_0^2 (t-t_0)^2} \exp & \left(- C^- \Psi(x,y+ x_0 \eps(t-t_0),t-\eps(t-t_0);x_0,y_0,t_0) \right) \le\\
   & \  \Gamma(x,y,t;x_0,y_0,t_0)  \le \\
   & \quad \qquad
   \frac{ C_{\eps}^+}{x_0^2 (t-t_0)^2} \exp \left(- c^+\Psi(x,y- x_0 \eps,t+ \eps;x_0,y_0,t_0) \right),
\end{split}
\end{equation}
for every $(x,y,t)\in \R^+ \times ]-\infty,y_0-x_0 \eps (t-t_0)[ \times ]t_0,T[$.
\end{theorem}

Note that, since the proof of Theorem \ref{th-main} is based on local estimates of the solution of $\L u=0$ and $\L$ is
\emph{locally} well approximated by the operator introduced in \eqref{Kolmogorov_op}, the diagonal bound in
\eqref{e-twosidedbounds1} agrees with the diagonal term of $\G_0$ in \eqref{Fund_sol}. Furthermore, the diagonal
estimate corresponds to the product of the standard deviations of the random variables $X_t$ and $Y_t$ defined in
\eqref{Yor process}. Indeed,
\begin{equation*}
\begin{split}
  \text{\rm Var} (X_t) & = x_0^2 \, e^{2t}\left(e^{2t} - 1 \right) = 2 x_0^2 \, t + o(t), \quad \text{as} \ t \to 0, \\
  \text{\rm Var} (Y_t) & = x_0^2 \left(\tfrac{1}{6} \left(e^{4t} - 1 \right) - \tfrac{2}{3} \left(e^{t} - 1
\right) -  \left(e^{t} - 1 \right)^2\right) = \tfrac{2}{3} \, x_0^2 \, t^3 + o(t^3), \quad \text{as} \ t \to 0.
\end{split}
\end{equation*}

\medskip

Clearly, the knowledge of the asymptotic behavior of the function $\Psi$ is crucial for the application of our Theorem
\ref{th-main}. In \cite{CPR}, it is shown that one can write the function $\Psi$ in terms of the function $g$ defined as
follows
\begin{equation*} 
  g(r)=\left\{
         \begin{array}{ll}
           \frac{\sinh(\sqrt{r})}{\sqrt{r}}, & \hbox{$r>0$,} \\
           1, & \hbox{$r=0$,} \\
           \frac{\sin(\sqrt{-r})}{\sqrt{-r}}, & \hbox{$-\pi^2<r<0$,}
         \end{array}
       \right.
\end{equation*}
and it is proven the following proposition
\begin{proposition} \label{prop-main}
For every $(x, y, t), (x_0, y_0, t_0) \in \R^+ \times \R^2$, with $t_0 < t$ and $y_0 > y$,
we have
\begin{equation*} 
\left\{
  \begin{array}{ll}
     &\Psi(x_1,y_1,t_1;x_0,y_0,t_0) = E(t_1-t_0)+\frac{4(x_1+x_0)}{y_0-y_1} - 4 \sqrt{E+\frac{4x_1x_0}{(y_0-y_1)^2}},\\
      & \hbox{if $E\ge -\frac{\pi^2}{t_1-t_0}$;} \\
     &\Psi(x_1,y_1,t_1;x_0,y_0,t_0) = E(t_1-t_0)+\frac{4(x_1+x_0)}{y_0-y_1} + 4 \sqrt{E+\frac{4x_1x_0}{(y_0-y_1)^2}},\\
     &\hbox{if $ -\frac{4 \pi^2}{t_1-t_0} <E <-\frac{\pi^2}{t_1-t_0}$.}
  \end{array}
\right.
\end{equation*}
where
\begin{equation*} 
  E = \frac{4}{(t-t_0)^2} g^{-1}\left(\frac{y_0-y}{(t - t_0) \sqrt{x x_0}}\right).
\end{equation*}
Moreover,
\begin{equation*} 
  \frac{\Psi(x,y,t;x_0,y_0,t_0)}{\frac{4}{(t - t_0)}\log^2\big(\tfrac{y_0 - y}{(t - t_0)
\sqrt{xx_0}}\big)+ \frac{4(x_0+x)}{y_0-y}} \to 1,
  \qquad \text{as} \qquad \frac{y_0-y}{(t - t_0)\sqrt{x_0x}} \to + \infty;
\end{equation*}
\begin{equation*} 
  \frac{\Psi(x,y,t;x_0,y_0,t_0)}{ \frac{4(\sqrt{x}+\sqrt{x_0})^2}{y_0-y}-\tfrac{ 4\pi^2}{(t - t_0)}} \to 1,
  \qquad \text{as} \qquad \frac{y_0-y}{(t - t_0)\sqrt{x_0x}} \to 0.
\end{equation*}
\end{proposition}

The above expression for the value function $\Psi$ has been obtained by using the
Pontryagin Maximum Principle \cite{PMP}, the upper bound in \eqref{e-twosidedbounds1} is a consequence of the fact that $\Psi$ satisfies the
Hamilton-Jacobi-Bellman equation $Y \Psi + \frac14 (X \Psi)^2 = 0$.

\medskip

To our knowledge, it is not easy to compare the integral expression of $p$ in \eqref{e-Yor-density} with the estimates
given in Proposition \ref{prop-main}, then Theorem \ref{th-main} provides us with an alternative explicit information
on the asymptotic behavior of $p$. Moreover, the method described in this section also applies to the divergence form
operator $\widetilde \L$ defined as
$$
  \widetilde \L u =x \p_{x} \left( a\, x \p_{x} u \right)  +b\, x\p_x u +x\p_y u - \p_t u,
$$
where $a$ and $b$ are smooth bounded coefficients, with $a$ bounded by below and $x \p_x a$ bounded. Note that, in this
case, an expression of $\Gamma$ analogous to \eqref{e-Yor-density} is not available. A further consequence of
\eqref{e-twosidedbounds1} is the following result. By applying \eqref{e-twosidedbounds1} to $\Gamma$ and to the
fundamental solutions $\Gamma^\pm$ of the operators
\begin{equation} \label{e-asian-pm}
  \L^{\pm} u =  \lambda^{\pm} x^2 \partial_{xx} u + x\p_x u+ x \partial_{y} u - \partial_{t} u,
  \quad (x,y,t) \in \R^+ \times \R \times ]0,T[,
\end{equation}
we obtain
\begin{align*} 
  {k^-}\Gamma^- &\big(x,y + \eps (t+1), t - \eps(t+1)\big)\\
& \le \Gamma(x,y,t) \\
& \quad \le {k^+}\Gamma^+ \Big(x,y -\tfrac{\eps}{1-\eps}(t+1), t+\tfrac{\eps }{1-\eps}(t+1)\Big),
\end{align*}
for every $(x,y,t), \in \R^+ \times \R \times ]0,T[$ with $y + \eps (t+1) < 0$ and $t>\eps /(1-\eps)$.
Hence, we obtain lower and upper bounds for the fundamental solution $\G$ of the \emph{variable coefficients operator}
$\widetilde \L$ in terms of the fundamental solutions $\G^{\pm}$ of the \emph{constant coefficients operators}
$\L^{\pm}$, whose expressions, up to some scaling parameters, agree with the function $p$ in \eqref{e-Yor-density}. We
refer to the article \cite{CPR} for the precise statement of the results of this section and for further details.

\subsubsection*{Acknowledgement}
{We thank the anonymous referee for his/her careful reading of our manuscript and for several suggestions
that have improved the exposition of our work.

\def\cprime{$'$} \def\cprime{$'$} \def\cprime{$'$}
  \def\lfhook#1{\setbox0=\hbox{#1}{\ooalign{\hidewidth
  \lower1.5ex\hbox{'}\hidewidth\crcr\unhbox0}}} \def\cprime{$'$}
  \def\cprime{$'$} \def\cprime{$'$}

\end{document}